\documentclass{article}
\usepackage[]{amsmath}
\usepackage[]{amsthm}
\usepackage[]{seqsplit}
\usepackage[]{tikz}
\usepackage[]{pgfplots}
\usepackage[]{xcolor}
\usepackage{xurl}
\usepackage[]{hyperref}

\pgfplotsset{compat=1.18}

\newtheorem{theorem}{Theorem}
\theoremstyle{definition}
\newtheorem{definition}{Definition}

\title{Cutting 4 by $n$ grids into two congruent pieces}
\author{Robert Dougherty-Bliss, Natalya Ter-Saakov, and  Doron Zeilberger}
\date{\today}

\def\1{{\overline{1}}}
\def\2{{\overline{2}}}
\def\frac#1#2{{#1 \over #2}}

\definecolor{pastelyellow}{HTML}{ffe45e}
\definecolor{pastelblue}{HTML}{7fc8f8}

\newcommand{\binarymatrix}[3]{%
  \begin{tikzpicture}[scale=0.5]
    \foreach \row [count=\y from 0] in {#3} {
      \foreach \val [count=\x from 0] in \row {
        \ifnum\val=1
          \fill[pastelyellow] (\x,-\y) rectangle ++(1,-1);
        \else
          \fill[pastelblue] (\x,-\y) rectangle ++(1,-1);
        \fi
      }
    }
    \draw[step=1,black,thick] (0,0) grid (#2,-#1);
  \end{tikzpicture}%
}

\begin{document}

\maketitle

{\it Pour Jean-Paul Delahaye, avec admiration}

\section*{Introduction} In the March 2025 issue of {\it Pour La Science} (the
French analogue of \emph{Scientific American}), in his delightful monthly
column [D], Jean-Paul Delahaye described in detail the solution of the
following counting problem.

{\it In how many ways can you cut a $3 \times 2n$ rectangle consisting of $6n$
    unit squares (in other words a $3 \times 2n$ (colorless) checkerboard) into
two \textbf{connected, congruent pieces?}}

In collaboration with his {\it \'epouse}, Martine Raison, Delahaye used human
ingenuity to prove that this number is given {\bf explicitly} by the beautiful
formula $2^{n+1}-n-1$. Here are all twelve ways to divide $3 \times 6$
rectangles into two congruent, connected parts:

\begin{center}
\binarymatrix{3}{6}{
{1, 1, 1, 1, 1, 1},
{0, 0, 0, 1, 1, 1},
{0, 0, 0, 0, 0, 0}}\quad
\binarymatrix{3}{6}{
{1, 1, 1, 1, 1, 1},
{1, 0, 0, 1, 1, 0},
{0, 0, 0, 0, 0, 0}}\quad
\binarymatrix{3}{6}{
{1, 1, 1, 1, 1, 1},
{0, 1, 0, 1, 0, 1},
{0, 0, 0, 0, 0, 0}}

\binarymatrix{3}{6}{
{1, 1, 1, 1, 1, 1},
{1, 1, 0, 1, 0, 0},
{0, 0, 0, 0, 0, 0}}\quad
\binarymatrix{3}{6}{
{1, 1, 1, 1, 1, 1},
{0, 0, 1, 0, 1, 1},
{0, 0, 0, 0, 0, 0}}\quad
\binarymatrix{3}{6}{
{1, 1, 1, 1, 1, 1},
{1, 0, 1, 0, 1, 0},
{0, 0, 0, 0, 0, 0}}

\binarymatrix{3}{6}{
{1, 1, 1, 1, 1, 1},
{0, 1, 1, 0, 0, 1},
{0, 0, 0, 0, 0, 0}}\quad
\binarymatrix{3}{6}{
{1, 1, 1, 1, 1, 1},
{1, 1, 1, 0, 0, 0},
{0, 0, 0, 0, 0, 0}}\quad
\binarymatrix{3}{6}{
{0, 1, 1, 1, 1, 1},
{0, 0, 0, 1, 1, 1},
{0, 0, 0, 0, 0, 1}}

\binarymatrix{3}{6}{
{0, 1, 1, 1, 1, 1},
{0, 1, 0, 1, 0, 1},
{0, 0, 0, 0, 0, 1}}\quad
\binarymatrix{3}{6}{
{0, 1, 1, 1, 1, 1},
{0, 0, 1, 0, 1, 1},
{0, 0, 0, 0, 0, 1}}\quad
\binarymatrix{3}{6}{
{0, 1, 1, 1, 1, 1},
{0, 1, 1, 0, 0, 1},
{0, 0, 0, 0, 0, 1}}
\end{center}

We asked ourselves, what about a $4 \times n$ checkerboard?

Here are twelve such ways to split a $4 \times 6$ rectangle into two congruent
parts:

\begin{center}
\binarymatrix{4}{6}{
{1, 1, 1, 1, 1, 1},
{0, 0, 0, 0, 0, 1},
{0, 1, 1, 1, 1, 1},
{0, 0, 0, 0, 0, 0}}\quad
\binarymatrix{4}{6}{
{1, 1, 1, 1, 1, 1},
{0, 1, 0, 0, 0, 1},
{0, 1, 1, 1, 0, 1},
{0, 0, 0, 0, 0, 0}}\quad
\binarymatrix{4}{6}{
{1, 1, 1, 1, 1, 1},
{1, 1, 0, 0, 0, 1},
{0, 1, 1, 1, 0, 0},
{0, 0, 0, 0, 0, 0}}

\binarymatrix{4}{6}{
{1, 1, 1, 1, 1, 1},
{0, 0, 1, 0, 0, 1},
{0, 1, 1, 0, 1, 1},
{0, 0, 0, 0, 0, 0}}\quad
\binarymatrix{4}{6}{
{1, 1, 1, 1, 1, 1},
{0, 1, 1, 0, 0, 1},
{0, 1, 1, 0, 0, 1},
{0, 0, 0, 0, 0, 0}}\quad
\binarymatrix{4}{6}{
{1, 1, 1, 1, 1, 1},
{1, 1, 1, 0, 0, 1},
{0, 1, 1, 0, 0, 0},
{0, 0, 0, 0, 0, 0}}

\binarymatrix{4}{6}{
{1, 1, 1, 1, 1, 1},
{0, 1, 0, 1, 0, 1},
{0, 1, 0, 1, 0, 1},
{0, 0, 0, 0, 0, 0}}\quad
\binarymatrix{4}{6}{
{1, 1, 1, 1, 1, 1},
{1, 1, 0, 1, 0, 1},
{0, 1, 0, 1, 0, 0},
{0, 0, 0, 0, 0, 0}}\quad
\binarymatrix{4}{6}{
{1, 1, 1, 1, 1, 1},
{0, 1, 1, 1, 0, 1},
{0, 1, 0, 0, 0, 1},
{0, 0, 0, 0, 0, 0}}

\binarymatrix{4}{6}{
{0, 0, 0, 1, 1, 1},
{0, 0, 0, 1, 1, 1},
{0, 0, 0, 1, 1, 1},
{0, 0, 0, 1, 1, 1}}\quad
\binarymatrix{4}{6}{
{0, 0, 0, 1, 1, 1},
{0, 0, 1, 1, 1, 1},
{0, 0, 0, 0, 1, 1},
{0, 0, 0, 1, 1, 1}}\quad
\binarymatrix{4}{6}{
{0, 0, 0, 1, 1, 1},
{0, 1, 1, 1, 1, 1},
{0, 0, 0, 0, 0, 1},
{0, 0, 0, 1, 1, 1}}
\end{center}

Now human ingenuity is still {\it necessary}, but (most probably) not {\it
sufficient}. In collaboration with our beloved silicon friends we proved the
following deep theorem.

\begin{theorem}
   Let $c_n$ be the number of ways of cutting a $4 \times n$ grid into two
   (connected) congruent pieces. The (ordinary) generating function of $c_n$
   is
$$
\sum_{n=0}^{\infty} c_n  x^n \, = \,
\frac{p(x)}{\left(1-2 x^{2}-x^{4}\right) \left(1-x \right)^{2} \left(1+x\right) \left(1-3 x^{2}-x^{4}\right)},
$$ where 
        \begin{align*}
            p(x) =
            &\ x(1+2x-4x^2-11x^3+11x^4+38x^5-35x^6-50x^7+50x^8\\
            &\quad -5x^9+5x^{10}+24x^{11}-24x^{12}+6x^{13}-6x^{14})
        \end{align*}
    The coefficients have the (approximate) asymptotic expansion
    \begin{equation*}
        c_n \sim 1.78631 (1 + 0.0842 (-1)^n) (1.8174)^n.
    \end{equation*}
\end{theorem}

For the sake of the OEIS the first $30$ terms are:
\begin{align*}
&1, 3, 5, 6, 20, 46, 76, 171, 277, 615, 981, 2158, \\
&3404, 7440, 11644, 25335, 39433, 85525, 132601, \\
&286944, 443664, 958524, 1479124, 3191905, 4918545, \\
&10605207, 16325329, 35178876, 54113204, 116555530
\end{align*}

To our delight, this sequence was not (as of Sept.\ 4, 2025) in the OEIS. We
hope to submit it soon.

Our colorful rectangles are very pretty, but it is more convenient to represent
such {\it decoupages} as $4 \times n$ matrices with entries $0$ and $1$, where
the value at position $ij$ indicates to which part that entry belongs. For
example:
\begin{equation*}
\vcenter{\hbox{
\binarymatrix{4}{6}{
{0 , 0 , 0 , 1 , 1 , 1},
{0 , 0 , 0 , 0 , 0 , 1},
{0 , 1 , 1 , 1 , 1 , 1},
{0 , 0 , 0 , 1 , 1 , 1}}}}
\to
\begin{bmatrix}
0 & 0 & 0 & 1 & 1 & 1 \\
 0 & 0 & 0 & 0 & 0 & 1 \\
  0 & 1 & 1 & 1 & 1 & 1 \\
   0 & 0 & 0 & 1 & 1 & 1
\end{bmatrix}
\end{equation*}

\begin{definition}
    A $m \times n$ matrix $M$ with entries in $\{0, 1\}$ that satisfies the
    rule
    \begin{equation}
        \label{rule}
        M_{ij} = 1 - M_{(m - i)(n - j)}
    \end{equation}
    such that the $0$'s and $1$'s form exactly two connected components (using
    north-south-east-west adjacency) is called a \emph{Graham matrix}. (Note
    that this rule uses \emph{zero} indexing.)
\end{definition}

Why Graham? For one, to honor the great Ron Graham, who advanced discrete
mathematics in a substantial way; for another, because the problem is something
like breaking apart a graham cracker fairly.

Note that if we drop the condition of the two pieces being congruent, and
consider only square matrices, then we have the (computationally) challenging
problem of the {\it Gerrymander} sequence [Sl], beautifully treated in [KKS]
and extended in [GJ]. 

\subsection*{How did we find this beautiful generating function?}

We asked our computer to show us all these creatures up to $n=12$. After some
thought, we convinced ourselves that Graham matrices could be recognized by
reading them one column at a time; in other words, that they could be
recognized by a \emph{finite state machine}, or a \emph{regular grammar}.

A finite state machine (also called a \emph{deterministic finite automata}) is
an object which detects words of a language. It reads in words one ``symbol''
at a time, and uses these symbols to determine a walk on a finite graph. Some
states in this graph are marked as \emph{accepting}, and all others are
\emph{rejecting}. A string of symbols is accepted by the machine if the walk it
determines ends in an accepted state.

\begin{figure}[t]
    \centering
    \includegraphics[width=\textwidth]{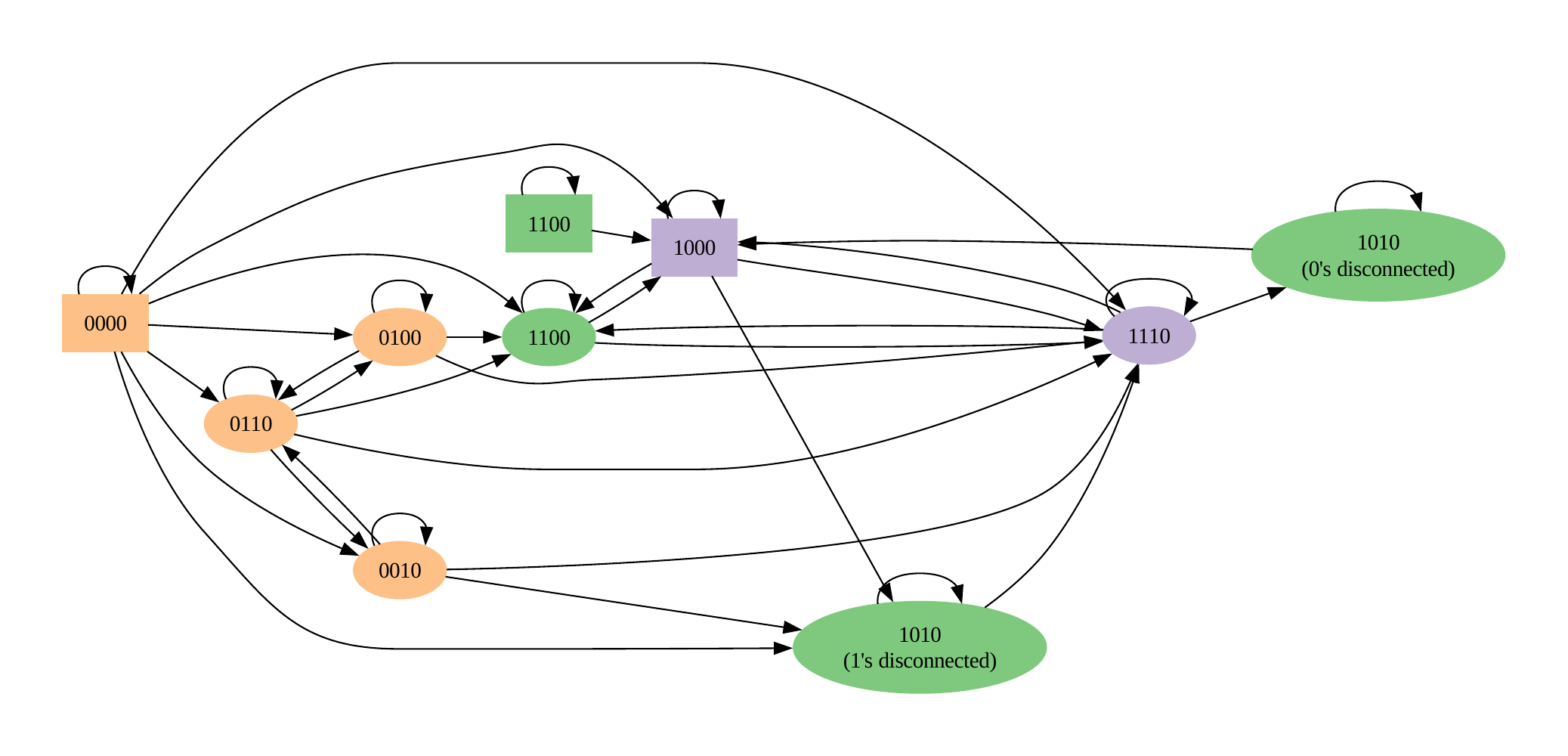}
    \includegraphics[width=0.7\textwidth]{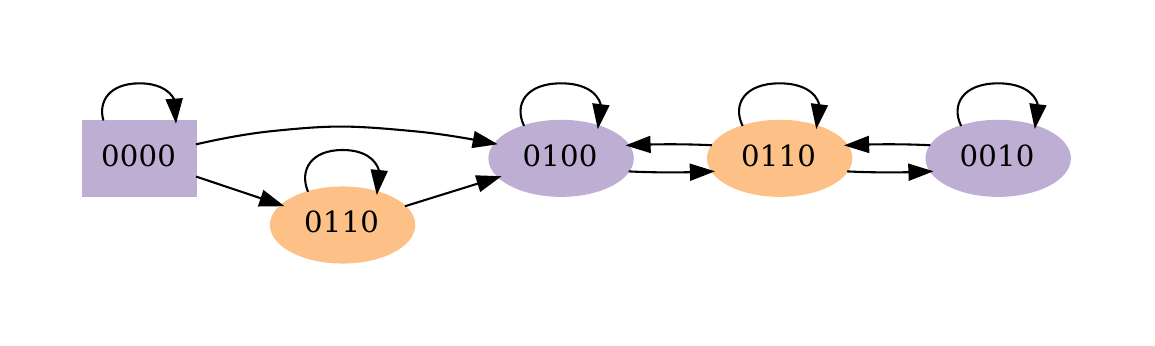}
    \includegraphics[width=0.5\textwidth]{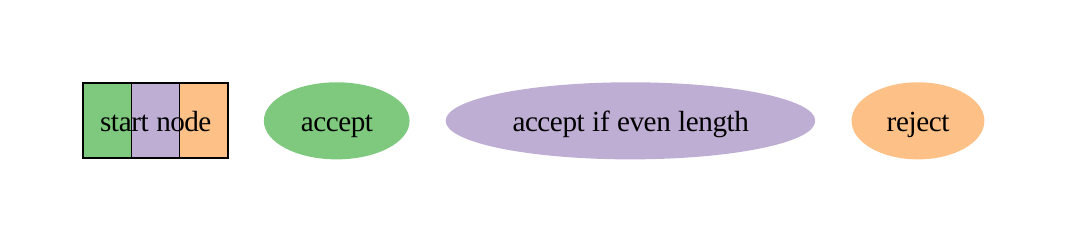}
    \caption{Finite state machine constructed to recognize $4 \times n$ Graham
    matrices. Rectangular states are start states, purple states accept if the
    string has an even number of symbols, and green states accept no matter
    what.}
    \label{fsm}
\end{figure}

In principle, it is not hard to see that Graham matrices can be recognized by a
finite state machine. We can track connectivity data and update it whenever a
column is read in; when the symbols stop, we apply rule~\eqref{rule} to
generate the second half of the matrix, merge the implied connectivity data,
then check that everything works. The difficult part is actually
\emph{constructing} this machine. In this problem it turns out that there are
shortcuts about the connectivity data which make the graph a more reasonable
size. If you are in a hurry, see Figure~\ref{fsm}.

The motivation behind constructing a finite state machine is the \emph{transfer
matrix method} ([St], section 4.7). If a finite state machine has adjacency
matrix $T$ (sometimes called a transfer matrix), then the \emph{symbolic
inverse} $(I - xT)^{-1}$ contains the generating functions which count walks
between different vertices of the state machine. The generating function we
want is
\begin{equation*}
    \sum_{\substack{\text{start state $i$} \\ \text{accepting state $j$}}} [(I - xT)^{-1}]_{ij}.
\end{equation*}

We will now actually describe the grammar.

\section*{The grammar}

A matrix that follows the complement rule \eqref{rule}
\begin{equation*}
    M_{ij} = 1 - M_{(4 - i)(n - j)}
\end{equation*}
is determined by its left half. To avoid double-counting symmetric
arrangements, we need to stipulate some conditions.
\begin{enumerate}
    \item The left-half of the bottom row of $M$ is all $0$'s and the first row to have a different number of $0$'s and $1$'s has more $0$'s.
    \item The first column of $M$ with a different number of $0$'s and $1$'s has at least as many $0$'s as $1$'s.
\end{enumerate}
We can enforce these conditions by a combination of rotations and reflections.
First, force a $0$ in the bottom left corner by rotating or reflecting. Then,
there are two cases:
\begin{enumerate}
    \item If the bottom-right corner is a $0$, then the entire bottom row must
        be $0$ to satisfy the connectivity rule. Then, either the first chunk or last chunk of
        columns of $M$ satisfy the second condition, and we can force it to be
        the first chunk with a reflection.

    \item If the bottom-right corner is a $1$, then the top-left corner is a
        $0$, and the first column is all $0$'s. The bottom and top rows consist
        of a run of $0$'s followed by a run of $1$'s. If the run of $0$'s in
        the bottom row is less than half of the row, then the run of $0$'s in
        the top row is at least half the row by rule~\eqref{rule}. Reflect if
        necessary to make this happen in the bottom row. If the bottom row is half and half, repeat this arguement with the first row that isn't.
\end{enumerate}

With these stipulations, there are only $2^3 = 8$ possible column vectors on
the left-hand side:
\begin{equation*}
\begin{bmatrix}
0 \\
0 \\
0 \\
0
\end{bmatrix}
\begin{bmatrix}
0 \\
0 \\
1 \\
0
\end{bmatrix}
\begin{bmatrix}
0 \\
1 \\
0 \\
0
\end{bmatrix}
\begin{bmatrix}
0 \\
1 \\
1 \\
0
\end{bmatrix}
\begin{bmatrix}
1 \\
0 \\
0 \\
0
\end{bmatrix}
\begin{bmatrix}
1 \\
0 \\
1 \\
0
\end{bmatrix}
\begin{bmatrix}
1 \\
1 \\
0 \\
0
\end{bmatrix}
\begin{bmatrix}
1 \\
1 \\
1 \\
0
\end{bmatrix}
\end{equation*}

According to the stipulations, the possible initial columns are $(0, 0, 0, 0)$,
$(1, 0, 0, 0)$, and $(1, 1, 0, 0)$. These columns all start with the 0's and
1's connected. The way we generate edges in Figure~\ref{fsm} is to look at all
possible next columns and track connectivity information. Almost all edges are
``trivial,'' meaning that the connectivity information does not change. For
example, here are four potential edges determined which begin at $(1, 1, 0,
0)$:
\begin{align*}
    &\begin{bmatrix}
        1 \\
        1 \\
        0 \\
        0
    \end{bmatrix} \to^?
    \begin{bmatrix}
        1 \\
        1 \\
        0 \\
        0
    \end{bmatrix} \text{--- yes; values remain connected} \\
    &\begin{bmatrix}
        1 \\
        1 \\
        0 \\
        0
    \end{bmatrix} \to^?
    \begin{bmatrix}
        0 \\
        0 \\
        1 \\
        0
    \end{bmatrix} \text{--- no; 1's are irreparably disconnected} \\
    &\begin{bmatrix}
        1 \\
        1 \\
        0 \\
        0
    \end{bmatrix} \to^?
    \begin{bmatrix}
        1 \\
        0 \\
        0 \\
        0
    \end{bmatrix} \text{--- yes; values remain connected} \\
    &\begin{bmatrix}
        1 \\
        1 \\
        0 \\
        0
    \end{bmatrix} \to^?
    \begin{bmatrix}
        1 \\
        1 \\
        1 \\
        0
    \end{bmatrix} \text{--- yes; values remain connected}
\end{align*}

The only time that connectivity information is changed is when we encounter the
column $(1, 0, 1, 0)$. \emph{Sometimes} when we see this column, new connected
components spawn. For example, in the edge
\begin{equation*}
    \begin{bmatrix}
        1 \\
        0 \\
        0 \\
        0
    \end{bmatrix} \to
    \begin{bmatrix}
        1 \\
        0 \\
        1 \\
        0
    \end{bmatrix},
\end{equation*}
we get a new $1$ connected component. Accordingly, there are two states in
Figure~\ref{fsm} for this column---one where the 1's are connected (and the 0's
disconnected), and one where the 1's are disconnected (and the 0's connected).

This explains all of the edges. But which states should be accepting?

To form a Graham matrix, the basic condition is that, when the second half of
the matrix is appended, the connected components need to ``line up.'' If an
even number of columns have been read in, then the next column of the final can
be determined from the current column by applying rule~\eqref{rule}. It is
routine to check whether the connected components line up, even in the
disconnected states, because the connectivity information also follows
rule~\eqref{rule}.

If an \emph{odd} number of columns have been read in, then the last column is
the ``middle'' column, and rule~\eqref{rule} must leave it fixed. The only
columns for which this is true are $(1, 1, 0, 0)$ and $(1, 0, 1, 0)$. In this
case, the first column in the right-half is determined by the
\emph{predecessor} of one of these columns. It is routine to check that no
matter what the predecessor is, the resulting matrix consists of exactly two
connected components.

We save the last word in this section for the lonely column $(0, 1, 1, 0)$, the
only column which is \emph{always} rejected.

\section*{Maple packages and asymptotics}

We have constructed our finite state machine in a Maple package {\tt
Decoupage.txt} that can be gotten from \url{https://sites.math.rutgers.edu/~zeilberg/tokhniot/Decoupage.txt}.

There are many ways to order the vertices of a finite state machine. We have
chosen, essentially, an arbitrary one. Our adjacency matrix is can be obtained
with \texttt{TransferMatrixMethods[transfer\_matrix]}. It is as follows:

\begin{equation*}
\left[\begin{array}{ccccccccccccccc}
1 & 1 & 1 & 1 & 1 & 1 & 1 & 0 & 0 & 1 & 0 & 0 & 0 & 0 & 0 
\\
 0 & 1 & 0 & 1 & 0 & 1 & 1 & 0 & 0 & 0 & 0 & 0 & 0 & 0 & 0 
\\
 0 & 0 & 1 & 1 & 0 & 1 & 0 & 0 & 0 & 1 & 0 & 0 & 0 & 0 & 0 
\\
 0 & 1 & 1 & 1 & 0 & 1 & 1 & 0 & 0 & 0 & 0 & 0 & 0 & 0 & 0 
\\
 0 & 0 & 0 & 0 & 1 & 1 & 1 & 0 & 0 & 1 & 0 & 0 & 0 & 0 & 0 
\\
 0 & 0 & 0 & 0 & 1 & 1 & 1 & 0 & 1 & 0 & 0 & 0 & 0 & 0 & 0 
\\
 0 & 0 & 0 & 0 & 1 & 1 & 1 & 0 & 0 & 0 & 0 & 0 & 0 & 0 & 0 
\\
 0 & 0 & 0 & 0 & 1 & 0 & 0 & 1 & 0 & 0 & 0 & 0 & 0 & 0 & 0 
\\
 0 & 0 & 0 & 0 & 1 & 0 & 0 & 0 & 1 & 0 & 0 & 0 & 0 & 0 & 0 
\\
 0 & 0 & 0 & 0 & 0 & 1 & 0 & 0 & 0 & 1 & 0 & 0 & 0 & 0 & 0 
\\
 0 & 0 & 0 & 0 & 0 & 0 & 0 & 0 & 0 & 0 & 1 & 0 & 1 & 0 & 0 
\\
 0 & 0 & 0 & 0 & 0 & 0 & 0 & 0 & 0 & 0 & 1 & 1 & 1 & 0 & 0 
\\
 0 & 0 & 0 & 0 & 0 & 0 & 0 & 0 & 0 & 0 & 0 & 0 & 1 & 1 & 0 
\\
 0 & 0 & 0 & 0 & 0 & 0 & 0 & 0 & 0 & 0 & 0 & 0 & 1 & 1 & 1 
\\
 0 & 0 & 0 & 0 & 0 & 0 & 0 & 0 & 0 & 0 & 0 & 0 & 0 & 1 & 1 
\end{array}\right]
\end{equation*}

The generating function matrix $(I - xM)^{-1}$ is too large to display here,
but the LCM of the denominators of its entries is
\begin{equation*}
\left(x^{2}-x +1\right) \left(x^{2}+3 x -1\right) \left(-1+x \right)^{3} \left(x^{2}+2 x -1\right).
\end{equation*}

Using Maple's \texttt{RootOf} representation in concert with partial fraction
decomposition, we can determine the asymptotic behavior of $c_n$ starting from
its generating function
\begin{equation*}
\frac{p(x)}{\left(-1+x \right)^{2} \left(x +1\right) \left(x^{4}+2 x^{2}-1\right) \left(x^{4}+3 x^{2}-1\right)},
\end{equation*}
where
\begin{align*}
    p(x) &=
    -6 x^{15}+6 x^{14}-24 x^{13}+24 x^{12}+5 x^{11}-5 x^{10}+50 x^{9} \\
    &\quad -50 x^{8}-35 x^{7}+38 x^{6}+11 x^{5}-11 x^{4}-4 x^{3}+2 x^{2}+x.
\end{align*}
The smallest poles of this rational function are $z$ and $-z$, where
\begin{equation*}
    z^{-1} = \frac{2}{\sqrt{2 \sqrt{13} - 6}} = 1.817354021\dots
\end{equation*}
The asymptotic formula is
\begin{align*}
    c_n &\sim \frac{1}{w} \left(\frac{22}{117} z^{2}+\frac{46}{117}+\frac{34}{117} z^{3}+\frac{103}{117} z\right) z^{-n} \\
    &\qquad + \frac{1}{w} \left(\frac{22}{117} z^{2}+\frac{46}{117}-\frac{34}{117} z^{3}-\frac{103}{117} z\right) (-z)^{-n} \\
        &\approx 1.7863 (1 + 0.08417 (-1)^n) z^{-n}.
\end{align*}

The front of this article

{\tt https://sites.math.rutgers.edu/\~{}zeilberg/mamarim/mamarimhtml/decoupage.html} contains some sample output files.

To see all $4 \times n$ Graham matrices for $1 \leq n \leq 12$ see the output
file:

{\tt https://sites.math.rutgers.edu/\~{}zeilberg/tokhniot/oDecoupage1.txt} \quad .

\bigskip

\section*{Future work}

What about the analogous problem of counting $5 \times 2n$, $6 \times n$, $7
\times 2n$, $8 \times n$ Graham matrices? As we mentioned earlier, it is not
hard to see that for any \emph{fixed} $m$, a finite state machine can recognize
the matrices of size $m \times n$. It is another matter to actually
\emph{construct} these machines.

We plan to follow-up with this question in a separate publication where
computers do everything---discover the grammar, process it, compute the
generating function, and so on. But even computers can only go so far, and we
are sure that humankind, and perhaps even computerkind will {\bf never} compute
the generating function for the sequence counting $100 \times n$ Graham
matrices.

{\bf{Acknowledgments.}} We would like to thank Gabriel Gendler for pointing out harder-to-notice symmetries. We would also like to thank Christoph Koutschan for helpful commentary and feedback.

{\bf References}

[D] Jean-Paul Delahaye, {\it Combinatoire pour les rectangles}, {\it Logique \& Calcul}, {\bf Pour La Science} No. {\bf 569}, Mars 2025.

[GJ] Anthony J. Guttmann and Iwan Jensen, {\it The gerrymander sequence, or A348456}, arXiv:2211.14482 [math.CO], 2022. {\tt https://arxiv.org/abs/2211.14482} \quad.

[KKS] Manuel Kauers, Christoph Koutschan, and George Spahn, {\it How Does the Gerrymander Sequence Continue?}, J. Int. Seq., {\bf 25} (2022), Article 22.9.7. \hfill\break
{\tt https://cs.uwaterloo.ca/journals/JIS/VOL25/Kauers/kauers6.html}.\hfill\break
arxiv version: {\tt https://arxiv.org/abs/2209.01787} \quad.

[Sl] Neil Sloane,  OEIS sequence A348456 {\tt https://oeis.org/A348456}.

[St] Richard P. Stanley,{\it ``Enumerative Combinatorics''}, Volume 1, Cambridge University Press, First Edition, Wadsworth and Brooks/Cole, 1986.

\bigskip
\hrule
\bigskip
Robert Dougherty-Bliss, Department of Mathematics, Dartmouth College, 29 N. Main Street, 6188 Kemeny Hall, Hanover NH 03755-3551 \quad .
Email: {\tt robert dot w dot bliss at gmail dot com}

\bigskip
Natalya Ter-Saakov, Department of Mathematics, Rutgers University (New Brunswick), Hill Center-Busch Campus, 110 Frelinghuysen
Rd., Piscataway, NJ 08854-8019, USA. \hfill\break
Email: {\tt  nt399 at rutgers dot edu }   \quad .

\bigskip
Doron Zeilberger, Department of Mathematics, Rutgers University (New Brunswick), Hill Center-Busch Campus, 110 Frelinghuysen
Rd., Piscataway, NJ 08854-8019, USA. \hfill\break
Email: {\tt DoronZeil at gmail  dot com}   \quad .
\end{document}